\documentclass[a4paper,12pt]{article}
\usepackage{amsmath}
\usepackage{amssymb}
\newcommand{\reals}{\mathbb {R}}
\newcommand{\naturals}{\mathbb{N}}
\newcommand{\complexes}{\mathbb{C}}
\newcommand{\sP}[1]{{\mathfrak {sp}}(#1)}
\newcommand{\End}[1]{{\mathrm {End}}(#1)}
\newcommand{\mapping}{\rightarrow}
\newcommand{\inner}[2]{\left< {#1} , {#2} \right>}
\newcommand{\tensor}{\otimes}
\newcommand{\iso}{\cong}
\newcommand{\trans}{{\pitchfork}}
\newcommand{\deriv}[2]{{\mathrm{D}}#1_{#2}}
\newcommand{\inv}[1]{{#1}^{-1}}
\newcommand{\tang}[2]{{\mathrm{T}}_{#2}{#1}}
\newcommand{\Tang}[1]{{\mathrm{T}}{#1}}

\newcommand{\sphere}[1]{{\mathrm{S}}^{#1}}

\newcommand{\Skew}[2]{{\Lambda^{#2}}(#1)}
\newcommand{\FORM}[2]{{\Omega^{#2}}(#1)}
\newcommand{\id}{\,\mathrm{l}\!\!\!\mathrm{1}}
\newcommand{\goth}[1]{\mathfrak #1}
 \newenvironment{prf}{{\bf Proof}\relax\\}{\hfill$\blacksquare$\par\bigskip}
 \newenvironment{prfof}[1]{{\bf Proof of {#1}}\relax\\}{\hfill$\blacksquare$\par\bigskip}
 \newtheorem{defn1}{Definition}[section]
 \newtheorem{lem1}[defn1]{Lemma}
 \newtheorem{prop1}[defn1]{Proposition}
 \newtheorem{thm1}[defn1]{Theorem}
 \newtheorem{cor1}[defn1]{Corollary}

 \newtheorem{claim}[defn1]{Claim}
 \newtheorem{rmk1}[defn1]{Remark}
 \newenvironment{rmk}{{\bf Remark}\relax\\}{\hfill\par}
 
\usepackage{eucal}
\usepackage{epsfig}
\newcommand{\quat}{\mathbb{H}}
\newcommand{\I}{{\,\widehat{\textbf{i}}\,}}
\newcommand{\J}{{\,\widehat{\textbf{j}}\,}}
\newcommand{\K}{{\,\widehat{\textbf{k}}\,}}
\newcommand{\eps}{\varepsilon}
\newcommand{\ol}[1]{\overline{#1}}
\newcommand{\ul}[1]{\underline{#1}}
\renewcommand{\d}{\mathrm{d}}
\newcommand{\hook}{
\begin{picture}(10,10)(-1,0)
\put(0,0){\line(1,0){7}}
\put(7,0){\line(0,1){7}}
\end{picture}}
\newcommand{\wdot}{{\,.\!\!\!\wedge\,}}
\newcommand{\OG}[2]{\Omega_G^{#2}\left({#1}\right)}
\newcommand{\dg}{{\d_\goth{g}}}

\renewcommand{\exp}[1]{e^{{#1}}}
\newcommand{\g}{\goth{g}}
\newcommand{\gs}{{\g^*}}
\newcommand{\Ig}{{\int_\g}}
\newcommand{\Igs}{{\int_\gs}}

\newcommand{\Ings}{\int_{N\times\gs}}

\newcommand{\ginner}[2]{\left\{#1,#2\right\}}
\title{Equivariant Integration formul{\ae} in HyperK\"ahler Geometry}
\author{Jonathan Munn}
\date{\today}
\begin{document}
\maketitle
\begin{abstract}
Lisa Jeffrey and Frances Kirwan developed an integration theory for symplectic
reductions. That is, given a symplectic manifold with symplectic group action,
they developed a way of pulling the integration of forms on the reduction
back to an integration of group-equivariant forms on the original space.\par
We seek an analogue of the symplectic integration formula as developed
by for the hyperK\"ahler case. This is almost
straightforward, but we have to overcome such obstacles as the lack of a
hyper-Darboux theorem and the lack of compactness in the case of hyperK\"ahler
reduction.
\end{abstract}
\tableofcontents
\section{Introduction}
While investigating the behaviour of instantons on $\sphere{4}$ in
\cite{MU}, we found that we need a result that links integrals
over a hyperK\"ahler reduction to an integral over the original
space which is, more often than not, a  much simpler space. To do
this we use the notion of equivariant cohomology, in particular
from the de Rham viewpoint, see \cite{BGV}.\par One of the major
results of this theory is the localisation formula.
\begin{thm1}\label{locform}
Let the torus $T$ act on the compact manifold $M$ with fixed set $M_0$, then
for any $\alpha\in H_T^\bullet(M)$
\[
\int_M\alpha=\int_{M_0}\frac{\iota^*\alpha}{e(\nu(M_0))}
\]
where $\nu(M_0)$ is the normal bundle of $M_0$ in $M$ and $e(\nu(M_0))$ its
Euler class on $M_0$
\end{thm1}
For the proof see \cite{S} pp24-25.\par Before we develop the
equivariant integration theory for hyperK\"ahler manifolds, we
need to motivate ourselves with the symplectic case which,
although well-established in \cite{JK} and \cite{GK}, we present
in some detail since the hyperK\"ahler case will follow in an
entirely similar manner. We hope that the reader will find this a
helpful comparison.
\section{Equivariant Cohomology of Symplectic Manifolds}
\subsection{Symplectic Group Actions}
Let $(M,\omega)$ be a compact symplectic manifold with a Hamiltonian action of
the compact Lie group $G$ which has bi-invariant metric $\{\cdot,\cdot\}$.
Let $\mu:M\mapping \goth{g}^*$ be the moment map of
this action, and assume that $0$ is a regular value of $\mu$ so that
$\inv{\mu}(0)$ is a manifold and hence so is the Marsden-Weinstein quotient
$M/\!\!/G=\inv{\mu}(0)/G$. Our goal is to relate integrals (i.e cohomology)
over $M/\!\!/G$, with the induced symplectic structure $\omega_0$, to integrals
over $M$.\par
Let $N=\inv{\mu}(0)$, then
we have the principal bundle
\[
\pi:N\mapping N/G=M/\!\!/G.
\]
 Let $\xi_1\ldots\xi_s$ be an orthonormal basis of $G$ and $\theta=\theta_i\xi_i$ be a connection on $N\mapping M/\!\!/G$, then
\[
\Omega=\bigwedge_{i=1}^s\theta_i
\]
is the induced volume form when restricting to the fibres of $N$, hence
\[
\pi_*\Omega=\mathrm{vol}\ G.
\]
The main result is
\begin{thm1}\label{equivintegral}
\begin{eqnarray*}
\int_{M/\!\!/G}\exp{i\omega_0}&=&\lim_{t\mapping\infty}\frac{(i)^{s^2}}{\left({{2\pi}}\right)^s\mathrm{vol} G}\Ig[\d y]\left(\exp{\left(\frac{-|y|^2}{4t}\right)}\int_M\exp{i\omega+i\ginner{\mu}{y}}\right)
\end{eqnarray*}
\end{thm1}
where at all times we are using a version of the Berezin integral formalism,
that the
integral of a $k$-form over a $l$-dimensional sub-manifold is zero if
$k\ne l$.
\subsection{The Kirwan Map}
Let $M$ be a compact symplectic manifold and $G$ a compact
connected Lie Group which acts in a Hamiltonian fashion on M, and
let $\mu$ be a moment map. Suppose that $0$ is a regular value of
$\mu$. In \cite{KI}, Frances Kirwan details a method by which we
can show that there is a surjective map
\[
\kappa:\mathrm{H}_G^\bullet(M;\complexes)\mapping\mathrm{H}^\bullet(M/\!\!/G;\complexes)
\]
which is given as $\kappa=\iota^*\inv{(p^*)}$ where $p:\inv{\mu}(0)\mapping M/\!\!/G$ is the quotient  map with
\[
p^*:\mathrm{H}^\bullet(M/\!\!/G)\mapping \mathrm{H}_G^\bullet(\inv{\mu}(0))
\]
being the isomorphism, since $G$ acts locally freely on $\inv{\mu}(0)$, and
$\iota:\inv{\mu}(0)\hookrightarrow M$ inclusion.
This is proved using the Morse theory of the function $f=|\mu|^2$ using the
Killing norm on $\gs$.\par
We should like to obtain a similar result for the boundaries of symplectic
manifolds, and later consider hyperK\"ahler manifolds and boundaries of hyperK\"ahler manifolds.\par
The function $f$ gives us a stratification of $M$, that is a finite collection
of subspaces $S_\alpha=f^{-1}\left((-\infty,\alpha]\right)$ for each critical
value $\alpha$, such that
\begin{enumerate}
\item $$M=\bigcup_{\alpha\in A}S_\alpha,$$
\item there is a partial ordering $>$ on $A$ such that
$$ \ol{S_\alpha}\subseteq\bigcup_{\gamma\ge\alpha}S_\gamma.$$
\end{enumerate}
Let $M_\alpha=\bigcup_{\alpha\ge\gamma}S_\gamma$.
We have the result
\begin{lem1}[\cite{KI}, Lemma 2.18 pp33-34]\label{surjhom}
Suppose $\{S_\gamma|\gamma\in A\}$ is a smooth $G$-invariant stratification of
$M$ such that for each $\alpha\in A$, the equivariant Euler class of the normal
bundle to $S_\alpha$ in $M$, is not a zero divisor in
$$\mathrm{H}_G^\bullet(S_\alpha;\complexes).$$ Then the inclusion
\[
\iota: M_\alpha\backslash S_\alpha\hookrightarrow M_\alpha
\]
 induces a surjection
\[
\iota^*:\mathrm{H}_G^\bullet(M_\alpha)\mapping\mathrm{H}_G^\bullet(M_\alpha\backslash S_\alpha).
\]
\end{lem1}
For the proof see \cite{KI}. For symplectic manifolds with a
Hamiltonian action of $\sphere{1}$, we know that not only is the
fixed set a submanifold, but its Euler class is not a zero
divisor! Kirwan also proves that the stratifications obtained by
minimally degenerate functions are just as good, and that moment
maps are minimally degenerate and the stratifications satisfy the
hypothesis in Theorem \ref{surjhom}. Prof. Kirwan has also proved
surjectivity in the hyperK\"ahler case using the norm square of
the hyperK\"ahler moment map $|\vec{\mu}|^2$, see \cite{KI1}.
\subsection{The Symplectic Structure near $\inv{\mu}(0)$}
Guillemin and Sternberg in \cite{GS} give a proof of the coisotropic embedding
theorem which allows us to describe the structure of $M$ in a neighbourhood of
$N$. We state it in the following form
\begin{thm1}[The Coisotropic Embedding Theorem]\label{Coisoemb}\relax
 (Guillemin and Sternberg \cite{GS} p315)
Given a symplectic manifold $(B,\omega_0)$ and a principal $G$-bundle
$\pi:P\mapping B$ then there is a unique%
\footnote{upto symplectomorphism, the uniqueness is a local quality}
symplectic
manifold $(M,\omega)$ into which $P$ embeds as a coisotropic submanifold and
the restriction of $\omega$ to $P$ is precisely $\pi^*\omega_0$.
\end{thm1}
The main corollary to this is the following
\begin{cor1}[\cite{JK}]\label{Normalsympform}
Given a compact symplectic manifold $(M,\omega)$ with a Hamiltonian $G$-action
and
moment map $\mu:M\mapping\gs$ with $0$ a regular value and $N=\inv{\mu}(0)$,
there is a neighbourhood $U$ of $N$ and a diffeomorphism
$\Phi:N\times\ol{B_{h}(0;\gs)}\mapping U$ such that
\[
\Phi^*\omega=\pi^*\omega_0+\d\tau
\]
where $\tau_{(p,\phi)}=\ginner{\phi}{\theta_p}$, where $\theta$ is a connection
on $\pi:N\mapping M/\!\!/G$.
\end{cor1}
Corollary \ref{Normalsympform} tells us that functions and forms supported
on a sufficiently small neighbourhood of $N$ can be replaced by functions and
forms supported on a sufficiently small neighbourhood of $N\times\{0\}$ in
$N\times\gs$. This will prove very useful when we come to consider equivariant
integration.
\subsection{Equivariant Integration}
Theorem \ref{equivintegral} was really introduced by Witten in \cite{WI} and
the technique that we use to prove this theorem is by Jeffrey and Kirwan in
\cite{JK} although they go further and produce some results proving a
localisation formula.
\par
We need to set up some theory of Gaussian integrals and Fourier analysis on Lie
algebras.
 Given any $p$-form $\alpha$
on a manifold, we write
\[
\exp{\alpha}=\sum_{i=0}^\infty \frac{\alpha^n}{n!}.
\]
We have a technical result which we will need.
\begin{prop1}\label{normalpowers}
Let $U$ and $V$ be $s$-dimensional manifolds.If  $\{a_j\}_{j=1}^s\subset\FORM{U}{1}$ and $\{b_j\}_{j=1}^s\subset\FORM{V}{1}$ then
\[
\exp{\left(i\sum_{j=0}^s a_j\wedge b_j\right)}_{[2s]}={i^{s^2}}[a]\wedge[b]
\]
where $[a]=\bigwedge_{i=1}^{s}a_i$.
\end{prop1}
%\begin{prf}
This is a simple calculation.\par
%\begin{eqnarray*}
%\exp{\left(i\sum_{j=0}^s a_j\wedge b_j\right)}_{[2s]}
%&=&\frac{1}{s!}\left(i\sum_{j=0}^{s}a_j\wedge b_j\right)^s\\
%&=& \frac{i^s}{s!}\left[s!\bigwedge_{j=1}^s a_j\wedge b_j\right]\\
%&=& {i^s}(-1)^{\frac{s}{2}(s-1)}[a]\wedge[b]\\
%&=& {i^{s^2}}[a]\wedge[b].
%\end{eqnarray*}
%\end{prf}
We have a general Gaussian integral formula, the proof of which is a standard argument using the method of ``Completing the square'':
\begin{equation}\label{gengauss}
\Ig[\d x]\exp{\left(-a|x|^2+b\inner{x}{y}\right)}=\left(\sqrt{\frac{\pi}{a}}\right)^s\exp{\left(\frac{b^2}{4a}|y|^2\right)}.
\end{equation}
Putting  $b=i$ and $a=1/4t$ into this, we have
\begin{equation}\label{gaussdelta}
\delta(x)=\lim_{t\mapping\infty}\left(\sqrt{\frac{t}{\pi}}\right)^s\exp{\left(-t|x|^2\right)}
=\lim_{t\mapping\infty}\left(\frac{1}{{2\pi}}\right)^s\Ig [\d y]
\exp{\left(\frac{-|y|^2}{4t}+i\inner{x}{y}\right)}.
\end{equation}

\begin{prfof}{\ref{equivintegral}}
(Following \cite{JK})
We put a connection on $N$ (or a metric on $M$) with connection form
$\theta$.
The generalised form
\begin{eqnarray*}
\lim_{t\mapping\infty}\int_\g\exp{\left(\frac{-|y|^2}{4t}\right)}
\exp{i\omega+i\ginner{\mu}{y}}
\end{eqnarray*}
on $M$ is supported on $N=\inv{\mu}(0)$ by (\ref{gaussdelta}),
hence by Corollary \ref{Normalsympform}
\begin{eqnarray*}
& &\lim_{t\mapping\infty}\Ig[\d y]\exp{\left(\frac{-|y|^2}{4t}\right)}\int_M\exp{i\omega+i\ginner{\mu}{y}}\\
&=&\lim_{t\mapping\infty}
\Ig[\d y]\exp{\left(\frac{-|y|^2}{4t}\right)}\Ings\exp{i\pi^*\omega_0+\ginner{\d z}{\theta} +
i\ginner{z}{\d\theta}+i\ginner{z}{y}}\\
&=&\lim_{t\mapping\infty}i^{s^2}\Ig[\d y]\exp{\left(\frac{-|y|^2}{4t}\right)}\Ings\exp{i\pi^*\omega_0 +
i\ginner{z}{\d\theta}+i\ginner{z}{y}}[\d z]\Omega\\
&=&\lim_{t\mapping\infty}i^{s^2}\Ig[\d y]\exp{\left(\frac{-|y|^2}{4t}\right)}\Igs[\d z]\exp{i\ginner{z}{y}}\int_N\exp{i\pi^*\omega_0+i\ginner{z}{\d\theta}}\Omega\\
&=&\lim_{t\mapping\infty}i^{s^2}(-1)^{s^2}\Igs\exp{i\ginner{z}{y}}[\d z]\Ig[\d y]\exp{\left(\frac{-|y|^2}{4t}\right)}\int_N\exp{i\pi^*\omega_0+i\ginner{z}{\d\theta}}\Omega\\
&=&i^{s^2}(-1)^{s^2}\left({{2\pi}}\right)^s\Igs[\d z]\delta(z)\int_N\exp{i\pi^*\omega_0+i\ginner{z}{\d\theta}}\Omega\\
&=&(-i)^{s^2}\left({{2\pi}}\right)^s\int_N\exp{i\pi^*\omega_0}\Omega\\
&=&(-i)^{s^2}{\left({{2\pi}}\right)^s}{\mathrm{vol} G}\int_{M/\!\!/G}\exp{i\omega_0}.
\end{eqnarray*}
\end{prfof}
We now obtain some useful formula for the integration of general forms on
$M/\!\!/G$.
\begin{thm1}\label{genequivintegral}
Let $\eta\in\OG{M}{\bullet}$ be equivariantly closed and have representative $\eta_0\in\FORM{M/\!\!/G}{\bullet}$.
If $M$ has no boundary then
\begin{eqnarray*}
\int_{M/\!\!/G}\exp{i\omega_0}\eta_0&=&\lim_{t\mapping\infty}\frac{(i)^{s^2}}{\left({{2\pi}}\right)^s\mathrm{vol} G}\Ig[\d y]\exp{\left(\frac{-|y|^2}{4t}\right)}\int_M\exp{i\omega+i\ginner{\mu}{y}}\eta(y).
\end{eqnarray*}
\end{thm1}
\begin{prf}
The proof of this fact is very similar to the proof of Theorem \ref{equivintegral}, but there are a few technicalities to overcome first.\par
Suppose that $\eta=\sum_I\eta_I \phi_I$ where $I$ is a multi-index, $\{\phi_i\}_{i=1}^s$ is a basis for $\gs$ and $\eta_I\in\FORM{M}{\bullet}^G$.
We have
\begin{eqnarray*}
& &\Ig[\d y]\exp{\left(\frac{-|y|^2}{4t}\right)}\int_M\exp{i\omega+i\ginner{\mu}{y}}\eta(y)\\
&=&\Ig[\d y]\exp{\left(\frac{-|y|^2}{4t}\right)}\int_M\exp{i\omega}\Igs[\d\phi]
\exp{i\ginner{\phi}{y}}\delta(\phi-\mu)\eta(y)\\
&=&\sum_I\Ig[\d y]\exp{\left(\frac{-|y|^2}{4t}\right)}\int_M\exp{i\omega}
\Igs[\d\phi]\exp{i\ginner{\phi}{y}}\delta(\phi-\mu)\eta_I\phi^I(y)\\
&=&(-1)^{s^2}\sum_I(-i)^{|I|}\int_M\exp{i\omega}\eta_I\Igs[\d\phi]\delta(\phi-\mu)\Ig[\d y]\exp{\left(\frac{-|y|^2}{4t}\right)}\exp{i\ginner{\phi}{y}}i^{|I|} y_I.
\end{eqnarray*}
Now, notice that the last integral is a derivative of the Fourier transform of
\[
y\mapsto\exp{\left(\frac{-|y|^2}{4t}\right)}
\]
which in the limit as $t\mapping\infty$ is the delta-function
$\phi\mapsto\delta(\phi)$. Hence the generalised function
\[
\delta(\phi-\mu)\Ig[\d y]\exp{\left(\frac{-|y|^2}{4t}\right)}\exp{i\ginner{\phi}{y}}i^{|I|} y_I
\]
is supported only on $N$, which means that by Corollary \ref{Normalsympform},
we can replace $M$ with $N\times\gs$ in the original integral.\par
Now
we have
\[
\iota^*\eta=\dg\beta+\pi^*\eta_0
\]
for some $\beta\in\OG{M}{\bullet}$, where $\iota:\inv{\mu}(0)\hookrightarrow M$
is the inclusion (see \cite{GK}).
Hence
\begin{eqnarray*}
\Ig[\d y]\exp{\left(\frac{-|y|^2}{4t}\right)}\int_M\exp{i\omega+i\ginner{\mu}{y}}\eta&=&\Ig[\d y]\exp{\left(\frac{-|y|^2}{4t}\right)}\int_{N\times\gs}\exp{i\omega+i\ginner{\mu}{y}}\dg\beta(y)\\
&+&\Ig[\d y]\exp{\left(\frac{-|y|^2}{4t}\right)}\int_{N\times\gs}\exp{i\omega+i\ginner{\mu}{y}}\pi^*\eta_0.
\end{eqnarray*}
Now let us examine the first term here. From the above argument
\begin{eqnarray*}
& &\Ig[\d y]\exp{\left(\frac{-|y|^2}{4t}\right)}\int_{N\times\gs}\exp{i\omega+i\ginner{\mu}{y}}\dg\beta(y)\\
&=&\Ig[\d y]\exp{\left(\frac{-|y|^2}{4t}\right)}\int_{N\times\gs}\dg\left(\exp{i\omega+i\mu(y)}\beta\right)(y)\\
& & \hbox{since }\omega+\mu(y)\hbox{ is an equivariantly closed form,}\\
&=&\Ig[\d y]\exp{\left(\frac{-|y|^2}{4t}\right)}\int_{N\times\gs}\d\left(\exp{i\omega+i\mu(y)}\beta\right)(y)\\
&-&\Ig[\d y]\exp{\left(\frac{-|y|^2}{4t}\right)}\int_{N\times\gs} X_y\hook\left(\exp{i\omega+i\mu(y)}\beta\right)(y)\\
&=& 0.
\end{eqnarray*}
The last term here vanishes since the integrand does not contain differential
forms of top degree, the first by Stokes' theorem.

The rest of the proof now follows that of Theorem
\ref{equivintegral}.
\end{prf}
\begin{rmk}
Using Theorem \ref{surjhom}, we find that for every $\eta_0\in\FORM{M/\!\!/G}{\bullet}$,
there is a $\eta\in\OG{M}{\bullet}$ such that
\[
\iota^*\eta=\pi^*\eta_0+\dg\beta,
\]
where $\iota:\inv{\mu}(0)\hookrightarrow M$ is the inclusion and $\pi:\inv{\mu}(0)\mapping
M/\!\!/G$ is the quotient map.
\end{rmk}
We will simplify matters by following Witten's notation, and write
\[
\oint_M\alpha(\xi)=\frac{1}{\mathrm{vol}\,(G)}\lim_{t\mapping\infty}
\int_\g\d\xi\exp{-\frac{|\xi|^2}{4t}}\int_M\alpha(\xi).
\]
In this notation we may rewrite the result of Theorem \ref{genequivintegral}
as
\begin{eqnarray*}
\int_{M/\!\!/G}\exp{i\omega_0}\eta_0&=&\frac{(i)^{s^2}}{\left({{2\pi}}\right)^s}\oint_M\exp{i\omega+i\ginner{\mu}{y}}\eta(y).
\end{eqnarray*}
\subsection{Localisation and Residues}
First we observe that in the case of a circle action, by Theorem \ref{genequivintegral} we have a Fourier transform
\begin{eqnarray*}
\int_{\mathcal{M}}\exp{i\omega_0}\eta_0&=&\frac{i}{\left({{2\pi}}\right)
}\oint_M\exp{i\omega+i\ginner{\mu}{y}}\eta(y).\\
&=&\frac{i}{\sqrt{2\pi}\,\mathrm{vol}\,\sphere{1}}\mathbf{F}\left(z\mapsto\int_M\exp{i\omega+i\ginner{\mu}{z}}\eta(z)\right)(0),
\end{eqnarray*}
where $\mathbf{F}$ denotes Fourier transform.\par
We also have a rather interesting fact that follows from Fourier theory and
Theorem \ref{genequivintegral}, namely,
\begin{thm1}\label{Restorthm}
For sufficiently small $\zeta\in\reals$
\begin{eqnarray*}
\int_{\mathcal{M}}\exp{i\omega_0}\eta_0
&=&\lim_{t\mapping\infty}\frac{(i)}
{\left({{2\pi}}\right)\mathrm{vol} \sphere{1}}
\int_\reals\d y\exp{\left(\frac{-|y|^2}{4t}\right)}
\int_M\exp{i\omega+i\ginner{\mu}{y}}\eta(y)\\
&=&\lim_{t\mapping\infty}\frac{(i)}
{\left({{2\pi}}\right)\mathrm{vol} \sphere{1}}
\int_\reals\d y
\exp{\left(\frac{-|y|^2}{4t}\right)}
\int_M\exp{i\omega+i\ginner{(\mu-\zeta)}{y}}\eta(y)\\
&=&\lim_{t\mapping\infty}\frac{(i)}
{\left({\sqrt{2\pi}}\right)\mathrm{vol} \sphere{1}}
\mathbf{F}\left(y\mapsto
\exp{\left(\frac{-|y|^2}{4t}\right)}
\int_M\exp{i\omega+i\ginner{\mu}{y}}\eta(y)\right)(\zeta).
\end{eqnarray*}
\end{thm1}
So this last result means that
\[
\frac{i}{\sqrt{2\pi}}\mathbf{F}\left(z\mapsto\int_M\exp{i\omega+i\ginner{\mu}{z}}\eta(z)\right)
\]
is smooth on a neighbourhood of $0$. But by the localisation
theorem (Theorem \ref{locform}), we know that if $\eta$ is
equivariantly closed
\[
\int_M\exp{i\omega+i\ginner{\mu}{z}}\eta(z)=\int_{M_0}\frac{\iota_{M_0}^*\left(\exp{i\omega+i\ginner{\mu}{z}}\eta(z)\right)}{e(z)}
\]
where $e(z)$ is the Euler class of the normal bundle of $M_0$ in $M$. Thus,
by Proposition 8.7 of \cite{JK}
\begin{eqnarray*}
& &\mathbf{F}\left(z\mapsto\int_M\exp{i\omega+i\ginner{\mu}{z}}\eta(z)\right)(0)\\
&=&\lim_{\eps\mapping 0_+}\frac{1}{\sqrt{2\pi}}\int_{\reals-i\xi}\mathbf{F}(\chi)(\eps\psi)
\mathbf{F^2}\left(z\mapsto\int_{M_0}\frac{\iota_{M_0}^*\exp{i\omega+i\ginner{\mu}{z}}\eta(z)}{e(z)}\right)(\psi)\d\psi\\
&=&\lim_{\eps\mapping 0_+}\frac{1}{(2\pi)^{\frac{3}{2}}}\int_{\reals-i\xi}\mathbf{F}(\chi)(\eps\psi)\int_{M_0}\frac{\iota_{M_0}^*\left(\exp{i\omega+i\ginner{\mu}
{-\psi}}\eta(-\psi)\right)}{e(-\psi)}\d\psi
\end{eqnarray*}
where $\chi$ is a smooth positive function on $\reals$ with compact support
and $\xi\in\reals$. Professors Jeffrey and Kirwan show that because
\[
\mathbf{F}\left(z\mapsto\int_M\exp{i\omega+i\ginner{\mu}{z}}\eta(z)\right)
\]
is smooth at 0, the latter integral is independent of $\chi$ and $\xi$.
Now we may bring to bear the standard theory of complex analysis. We know that
$\mathbf{F}(\chi)$ is entire when extended to $\complexes$ because $\chi$ is
smooth and compactly supported on $\reals$ (Proposition 8.4 of \cite{JK}), we
also know that $\mu$ is constant and non-zero on $M_0$ since $0$ is a regular
value of $\mu$.
Hence using the Residue theorem from complex analysis
\begin{eqnarray*}
& &\lim_{\eps\mapping 0_+}\frac{1}{(2\pi)^{\frac{3}{2}}}\int_{\reals-i\xi}\mathbf{F}(\chi)(\eps\psi)\int_{M_0}\frac{\iota_{M_0}^*\left(\exp{i\omega+i\ginner{\mu}
{-\psi}}\eta(-\psi)\right)}{e(-\psi)}\d\psi\\
&=&-\frac{i}{\sqrt{2\pi}}\mathrm{Coeff}_{y^{-1}}\int_{M_0^+}\frac{\iota_{M_0^+}^*\exp{i\omega+i\ginner{\mu}
{y}}\eta(y)}{e(y)}
\end{eqnarray*}
where $M_0^+$ is the set of fixed points on which $\mu$ is
positive, and we always take the Laurent expansion of the
integrand at $y=0$. Thus we end up with the following theorem.
\begin{thm1}\label{Resthm0}
Let $M$ be compact symplectic manifold  with Hamiltonian action of
the circle such that $0$ is a regular value of the moment map
$\mu$. Let $\iota:\inv{\mu}(0)\hookrightarrow M$ be inclusion. If
$\eta\in\Omega_{\sphere{1}}(M)$ is equivariantly closed and
satisfies
\[
\iota^*\eta(y)=\pi^*\eta_0+\dg\beta(y)
\]
for some closed $\eta_0\in\Omega(M/\!\!/\sphere{1})$ and $\beta\in\Omega_{\sphere{1}}(M)$, then
\begin{eqnarray*}
\int_{\mathcal{M}}\exp{i\omega_0}\eta_0&=&\frac{i}{\left({2\pi}\right)
}\oint_M\exp{i\omega+i\ginner{\mu}{y}}\eta(y)\\
&=&\frac{1}{2\pi\mathrm{vol}\,\sphere{1}}\mathrm{Coeff}_{y^{-1}}\int_{M_0^+}\frac{\iota_{M_0^+}^*\left(\exp{i\omega+i\ginner{\mu}{y}}{\eta(y)}\right)}{e(y)}
\end{eqnarray*}
where we always take the Laurent expansion of the integrand at $y=0$.
\end{thm1}
\begin{cor1}\label{Restorthm0}
Let $M$ be compact symplectic manifold  with Hamiltonian action of
the torus $\mathbb{T}^k$ such that
$0\in\mathrm{Lie}(\mathbb{T}^k)$ is a regular value of the moment
map $\mu$. Choose a splitting of $\mathbb{T}^k$ into the product
of circles Let $\iota:\inv{\mu}(0)\hookrightarrow M$ be inclusion.
If $\eta\in\Omega_{\mathbb{T}^k}(M)$ is equivariantly closed and
satisfies
\[
\iota^*\eta(y_1,\ldots,y_k)=\pi^*\eta_0+\dg\beta(y_1,\ldots,y_k)
\]
for some closed $\eta_0\in\Omega(M/\!\!/\sphere{1})$ and $\beta\in\Omega_{\sphere{1}}(M)$, then
\begin{eqnarray*}
\int_{\mathcal{M}}\exp{i\omega_0}\eta_0&=&\left(\frac{i}{2\pi}\right)^k
\oint_M\exp{i\omega+i\ginner{\mu}{\ul{y}}}\eta(\ul{y})\\
&=&\left(\frac{1}{2\pi}\right)^k\frac{1}{\mathrm{vol}\,\mathbb{T}^k}\mathrm{Coeff}_{y_1^{-1},
\ldots,y_k^{-1}}\int_{M_0^+}\frac{\iota_{M_0^+}^*\left(\exp{i\omega+i\ginner{\mu}{\ul{y}}}{\eta(\ul{y})}\right)}{e(\ul{y})}
\end{eqnarray*}
where $M_0^+$ is the set of fixed points on which each component of the moment
map $\mu$ induced by the splitting is positive. We always take the Laurent expansion of the integrand at $\ul{y}=0$.
\end{cor1}
\section{Equivariant Integrals in HyperK\"ahler Reduction}
\subsection{Preliminaries}
Recall that a hyperK\"ahler manifold $(M,\vec{\omega})$ is a $4k$-dimensional
Riemannian manifold with three symplectic forms which form the components of the
$\sP{1}$-valued 2-form
\[
\vec\omega=\omega_\I\tensor\I+\omega_\J\tensor\J+\omega_\K\tensor\K
\]
 whose associated complex structures (resp $(\I,\J,\K)$) obey
$$\I^2=\J^2=\K^2=\I\J\K=-\id.$$ Let $M$ be compact and $G$ act on $M$ with a
tri-Hamiltonian
action, that is, the action is Hamiltonian with respect to each of the
symplectic forms. Then we have a threefold moment map
\[
\vec\mu=(\mu_\I,\mu_\J,\mu_\K):M\mapping\g\tensor_\reals\Im\quat
\]
such that
\[
\ginner{\d\vec{\mu}(X)}{\xi}=\vec\omega(X_\xi,X)
\]
for each $\xi\in\g, X\in\Gamma\left(\Tang{M}\right)$. For an element
$\vec{a}=a_0+a_1\I+a_2\J+a_3\K\in\quat$ and $X\in\tang{M}{p}$, we set
\[
\vec{a}X=a_0X+a_1\I X+a_2\J X+a_3\K X.
\]
If $0$ is a regular value of $\vec\mu$, we can form the hyperK\"ahler reduction
\[
\mathcal{M}=\inv{\vec\mu}(0)/G=M/\!\!/\!\!/\!\!/G.
\]
Now we can regard this reduction in several ways:
\begin{enumerate}
\item as the quotient of $N_0=\inv{\vec\mu}(0)$ by $G$,
\item if we fix a complex structure $\I$ then set $\omega_\complexes=\omega_\J+i
\omega_\K$ and $\mu_\complexes=\mu_\J+i\mu_\K$, then we can view
$M/\!\!/\!\!/\!\!/G$ as $N_\complexes/\!\!/G$ with respect to the K\"ahler form
$\omega_\I$ where $N_\complexes=\inv{\mu_\complexes}(0)$,
\item fix a complex structure $\I$ and take the GIT symplectic reduction of
$M$ by the complex group $G_\complexes$.
\end{enumerate}

Essentially we'd like to use the first point of view. It is independent of any
choices and the various spaces and groups involved will be compact, a fact that
the other choices do not share.\par
We do have \emph{not} have a hyperK\"ahler version of the Coisotropic embedding theorem due to lack of a hyper-Darboux theorem. However, we can make some
attempt at looking at the hyperK\"ahler structure near $\inv{{\vec{\mu}}}(0)$
in the simpler case of a torus action.
\ref{Normalsympform}
\begin{lem1}\label{HNormalsympform}
Let $(M,\inner{\cdot}{\cdot},\vec\omega)$ be  a hyperK\"ahler manifold  with a
tri-Hamiltonian
action of the $k$-dimensional torus $\mathbb{T}^k$, and moment map $\vec\mu:M\mapping\gs\tensor\sP{1}\iso\goth{t}^k\tensor\Im\quat$.
 Suppose that $0$ is a regular value of $\vec\mu$, so that $N=\inv{\mu}(0)$ is
a submanifold. Then there are $k$ differential 1-forms
$\phi^\alpha\in\FORM{M;\quat}{1}$ such that on a sufficiently small tubular
neighbourhood of $N$, we have
\[
\vec{\omega}=\vec{\omega}'-\phi^\alpha\wedge\ol{{\phi^{\alpha}}}
\]
where $\vec{\omega}'$ is closed and restricts on $N$ to $\pi^*\vec{\omega}_0$, the pullback
under the projection of the hyperK\"ahler structure $\vec{\omega}_0$ of
$M/\!\!/\!\!/\!\!/\mathbb{T}^k$. As a result $\phi^\alpha\wedge\ol{{\phi^{\alpha}}}$ is also closed.

\end{lem1}
\begin{prf}
Note that in a small tubular neighbourhood $U$ around $N$ in $M$,
$\deriv{\vec{\mu}}{p}$ is surjective for all $p\in U$. So restricted to $U$
\[
\Tang{M}\iso\ker\d\vec\mu\oplus(\ker\d\vec{\mu})^\perp.
\]
Now, examining $\ker\d\vec{\mu}$, we notice that, because the group acting is
a torus, we have for all $p\in U$, $\xi,\eta\in\goth{t}^k$
\begin{equation}\label{Lag}
\begin{array}{rcl}
\vec{\omega}\left(X_\xi(p),X_\eta(p)\right)&=&\{\vec{\mu}(p),[\xi,\eta]\}\\
                                           &=&0
\end{array}
\end{equation}
where $\{\cdot,\cdot\}$ denotes the invariant inner product on
$\goth{g}$. Denote the sub-bundle of $\ker\d\vec{\mu}$ by $V$
which form the vertical vector fields over $M/\mathbb{T}^k$, and set
$$H=\ker\d\vec{\mu}\cap (V)^\perp.$$ Now, let
$\{X_\alpha\}_{\alpha=1}^k$ be an orthonormal frame of $V$ over a
(possibly smaller) $U$. Then
\[
\Tang{M}=H\oplus V\oplus V_{\Im\quat}
\]
where $V_{\Im\quat}=\hbox{Span}\{qX_\alpha\vert q\in\Im\quat\}_{\alpha=1}^k$.
Now we can set
\begin{eqnarray*}
\phi_0^\alpha&=&\frac{1}{\sqrt{2}}\inner{X_\alpha}{\cdot}\\
\phi_1^\alpha&=&\frac{1}{\sqrt{2}}\inner{\I X_\alpha}{\cdot}\\
\phi_2^\alpha&=&\frac{1}{\sqrt{2}}\inner{\J X_\alpha}{\cdot}\\
\phi_3^\alpha&=&\frac{1}{\sqrt{2}}\inner{\K X_\alpha}{\cdot}\\
\end{eqnarray*}
and
\[
\phi^\alpha=\phi_0^\alpha+\I\phi_1^\alpha+\J\phi_2^\alpha+\K\phi_3^\alpha.
\]
Now setting $q_0=1,q_1=\I,q_2=\J,q_3=\K$, we have for all $\alpha,\beta$ and
$i,j\in\{0,1,2,3\}$
\begin{eqnarray*}
\inner{q_iX_\beta}{q_jX_\gamma}&=&\inner{q_j^*q_iX_\beta}{X_\gamma}\\
                            &=&\left\{\begin{array}{cc}
\omega_{q_j^*q_i}(X_\beta,X_\gamma) & i\neq j\\
\inner{X_\beta}{X_\gamma} &i=j
\end{array}\right.\\
&=&\delta_{ij}\delta_{\beta\gamma},
\end{eqnarray*}
since $\vec{\omega}$ vanishes on $V$ by (\ref{Lag}). Hence the vector fields
${\{q_iX_\alpha\}_{\alpha=1}^k}_{i=0}^{3}$ form an orthonormal frame for $H^\perp$ on $U$.
Using this we can show that
\[
\vec{\omega}(q_jX_\beta,q_kX_\gamma)=\left(\eps_{ijk}+\delta_{0j}\delta_{ik}-\delta_{0k}\delta_{jk}\right)\delta_{\beta\gamma}q_i
\]
where
\[
\eps_{ijk}=\left\{
\begin{array}{cc}
1 &(i,j,k) \hbox{ is an even permutation of } (1,2,3)\\
-1 &(i,j,k) \hbox{ is an odd permutation of } (1,2,3)\\
0 &\hbox{otherwise.}
\end{array}
\right.
\]
It can be readily checked that
\[
\phi^\alpha\wedge\ol{\phi^\alpha}(q_jX_\beta,q_kX_\gamma)=-\left(\eps_{ijk}+\delta_{0j}\delta_{ik}-\delta_{0k}\delta_{jk}\right)\delta_{\beta\gamma}q_i.
\]
Now let $\vec{\omega}'=\vec{\omega}+\phi^\alpha\wedge\ol{\phi^\alpha}$.
Now we know by construction that for $v\in H, X\in H^\perp$
we have $\vec{\omega}(v,X)=0$, since $qX\in H^\perp$ for all $q\in\quat$.
So $\vec{\omega}'$ is a purely horizontal form. By construction,
$\phi^\alpha\wedge\ol{\phi^\alpha}$ is purely vertical/normal.
Since
\[
0=\d\vec{\omega}=\d\vec{\omega}'-\d(\phi^\alpha\wedge\ol{\phi^\alpha})
\]
we know that
$\d\vec{\omega}'$ can have at most one vertical/normal component,
and $\d(\phi^\alpha\wedge\ol{\phi^\alpha})$ at most one horizontal component.
Hence, by comparing components, we have $\d\vec{\omega}'=0$ and
$\d(\phi^\alpha\wedge\ol{\phi^\alpha})=0$.

\end{prf}
\begin{rmk} Notice that the $\phi_0^\alpha$ gives rise to the connection on
$N\mapping M/\!\!/\!\!/\!\!/\mathbb{T}^k$ generated by the metric when we restrict to $N$.
\end{rmk}

Now, a hyperK\"ahler manifold has a canonical $4$-form, namely%
\footnote{We use the dot $.$ to denote the scalar product of quaternions following the
motivation from regarding $\Im\quat$ as $\reals^3$. Hence $\wdot$ is an
operation
\[
\wdot:\Skew{V}{\bullet}\tensor\quat\times\Skew{V}{\bullet}\tensor\quat\mapping\Skew{V}{\bullet}
\]
for any vector space $V$, given by combining the quaternionic dot product with
the wedge product.}
\[
\Omega=\vec{\omega}\wdot\vec{\omega}=\omega_\I\wedge\omega_\I+\omega_\J\wedge\omega_\J+
\omega_\K\wedge\omega_\K.
\]
Our aim is to mimic the construction of the symplectic version of Witten's
equivariant integral, but there are hidden dangers in this approach. If we try
the approach of taking a complex structure, then not only do we lose the
invariance under change of complex structure but we have problems localising
the Poincar\'e dual of $N_\complexes
=\mu_\complexes^{-1}(0)$. Other methods involve reducing by a non-compact group
which may be fine in the algebraic-geometric sense but the integration is not
very well-defined.
\begin{defn1}
Let  $M$ be a hyperK\"ahler manifold which possesses a\newline
tri-Hamiltonian action of the compact Lie group $G$ and
$\alpha\in\OG{M}{\bullet}$ be an equivariant form. We shall say
that $\alpha$ is associated with
$\alpha_0\in\FORM{\mathcal{M}}{\bullet}$ if
\[
\iota^*\alpha=\pi^*\alpha_0+\dg\beta
\]
where $\iota:\vec{\mu}^{-1}(0)\mapping M$ is inclusion and
$\pi:\vec{\mu}^{-1}(0)\mapping\mathcal{M}$ is the quotient map and
$\beta\in\OG{\vec{\mu}^{-1}(0)}{\bullet}$. In the case that
$\alpha$ and $\alpha_0$ are compactly supported, we shall say that
$\alpha$ is compactly associated with $\alpha_0$ if $\alpha$ is
associated with $\alpha_0$ as above, and the form $\beta$ is also
compactly supported.
\end{defn1}
As always, in the de Rham model, equivariant forms are regarded as
polynomials functions on the Lie algebra $\goth{g}$ with
$\FORM{M}{\mathrm{even}}$-coefficients.
\begin{prop1}
The equivariant form
\[
\Omega(y)=(\vec{\omega}+\ginner{\vec{\mu}}{y})\wdot(\vec{\omega}+\ginner{\vec{\mu}}{y})
\] is an equivariantly closed $\reals$-valued form.
\end{prop1}
\begin{prf}
For each $q\in\{\I,\J,\K\}$, we have
\[
\dg(\omega_q+\{\mu_q,\xi\})(\xi)=\{\d\mu_q,\xi\}-X_\xi\hook\omega_q=0.
\]
Hence
\[
y\mapsto(\vec{\omega}+\{\vec{\mu},y\})\wdot(\vec{\omega}+\{\vec{\mu},y\})=\sum_{q\in\{\I,\J,\K\}}(\omega_q+\{\mu_q,y\})\wedge(\omega_q+\{\mu_q,y\})
\]
is an equivariantly closed $\reals$-valued form.
\end{prf}
Now, there is an issue of compactness arising here. There are no
compact hyperK\"ahler manifolds which admit tri-hamiltonian
actions of Lie groups. However, we can talk about certain
submanifolds of non-compact hyperK\"ahler manifolds which reduce
under the action of a tri-hamiltonian group action to submanifolds
of the quotient.
\begin{defn1}
We shall say that the tri-Hamiltonian action of a Lie group $G$ on
a hyperK\"ahler manifold $M$ with respects a $G$-invariant
submanifold $L$ of $M$ if for each $p\in L$ any vector normal to
$L$ at $p$ lies in $\ker\deriv{\vec{\mu}}{p}$.
\end{defn1}
\begin{rmk}
For a group action to respect a submanifold $L$ of a hyperK\"ahler
manifold, it is necessary that the codimension of $L$ in $M$
should not exceed three times the dimension of the group.
\end{rmk}
\begin{thm1}
Let $M$ be a hyperK\"ahler manifold that admits a tri-hamiltonian
action of the compact Lie group $G$ with moment map $\vec{\mu}$.
Let $0\in\goth{g}\tensor\Im\quat$ be a regular value and suppose
further that the action of $G$ respects a $G$-invariant
submanifold $L$ , then
\[
\frac{L\cap\inv{{\vec{\mu}}}(0)}{G}
\]
is a submanifold of $M/\!\!/\!\!/\!\!/G$ with the same codimension
as $L$ in $M$ .
\end{thm1}
\begin{prf}
This is an exercise in tranversality. For each $p\in L
\cap\inv{{\vec{\mu}}}(0)$
\begin{eqnarray*}
\tang{L}{p}+\tang{\inv{{\vec{\mu}}}(0)}{p}&=&\tang{L}{p}+\ker\deriv{\mu}{p}.
\end{eqnarray*}
Since the action respects $L$, $\ker\deriv{\mu}{p}$ contains the
normal vectors at $p$, proving that $L\trans\inv{{\vec{\mu}}}(0)$.
Hence this is a $G$ invariant manifold of dimension $n-3\dim G$
and the result follows.
\end{prf}
We will prove a Witten-style formula for submanifolds of
hyperK\"ahler spaces that are respected by the group action. To do
this, we need a number of results.
\begin{prop1}\label{suptsq}
 For $n\in\naturals$ and $x\in\reals\backslash\{0\}$
\[
\lim_{t\mapping\infty}\int_\reals\exp{-\frac{y^2}{4t}+ixy}y^{\frac{n}{2}}\d y=0
\]
\end{prop1}
\begin{prf}
A calculation.
\begin{eqnarray*}
\int_\reals\exp{-\frac{y^2}{4t}+ixy}y^{\frac{n}{2}}&=&\exp{-tx^2}\int_\reals\exp{-\left(\frac{y}{2\sqrt{t}}-ix\sqrt{t}\right)^2}y^{\frac{n}{2}}\d y\\
&=&2\exp{-tx^2}\sqrt{t}\int_{\reals-ix\sqrt{t}}\exp{-z^2}(2\sqrt{t}z+2ixt)^{\frac{n}{2}}\d z\\
&=&2\exp{-tx^2}t^{\frac{n+1}{2}}\int_{\reals-i\mathrm{sign}(x)}\exp{-z^2}\left(2\frac{z}{\sqrt{t}}+2ix\right)^{\frac{n}{2}}\d z\\
& &\hbox{since either branch of the square root}\\
& &\hbox{is holomorphic away from 0}\\
&\mapping& 0
\end{eqnarray*}
as $t\mapping\infty$, provided $x\neq 0$.
\end{prf}
\begin{thm1}\label{hypintform}
Let  $(M,\vec{\omega})$ be a hyperK\"ahler manifold with a
tri-hamiltonian%
\footnote{therefore implying $M$ is not compact}
 action of $\sphere{1}$ and associated hyperK\"ahler
moment map $\vec{\mu}:M\mapping\Im\quat$. Suppose that $0$ is a
regular value for $\vec{\mu}$ and let
$\mathcal{M}=M/\!\!/\!\!/\!\!/\sphere{1}$ be the hyperK\"ahler
reduction of $M$ with hyperK\"ahler structure $\vec{\omega}_0$. If
$\eta_0\in
\FORM{\mathcal{M}}{\bullet} $ is compactly supported and
compactly associated with $\eta\in
\Omega_{\sphere{1}}^\bullet(M)$, then
we have
\begin{eqnarray*}
& & (\dim\mathcal{M}-\mathrm{deg}\eta_0+1)\int_\mathcal{M}\exp{i\vec{\omega}_0\wdot\vec{\omega}_0}
\eta_0\\
&=&\frac{1}{6\pi i\sqrt{2}}\oint_M
\exp{i\vec{\omega}\wdot\vec{\omega}+2i{\sqrt{y}}{\vec{\mu}}.\vec{\omega}+i
|{\vec{\mu}}|^2y}\eta(\sqrt{y}).
\end{eqnarray*}
\end{thm1}

\begin{prf}
Notice that
\[
\exp{2i{\sqrt{y}}{\vec{\mu}}.\vec{\omega}}\eta(\sqrt{y})
\]
is a polynomial in $\sqrt{y}$ due to the presence of a (usual)
form in the exponent, so we can apply Proposition \ref{suptsq} to
show that the integrand is supported on $N$. \par Also, by Lemma
\ref{HNormalsympform}, we have in a neighbourhood of $N$
\begin{eqnarray*}
\vec{\omega}\wdot\vec{\omega}&=&{\vec{\omega}'}\wdot{\vec{\omega}'}\\
&-& \vec{\omega}'\wdot\left(\phi\wedge\ol{\phi}\right)\\
&+& 6\phi_0\wedge\phi_1\wedge\phi_2\wedge \phi_3\\
\end{eqnarray*}
and by the remark under Lemma \ref{HNormalsympform},
 $\phi_0$ is a connection on $N\mapping\mathcal{M}$ modified to take real
values and scaled by a factor of $\frac{1}{\sqrt{2}}$. From this
it is worth trying to estimate
$
\exp{i\vec{\omega}\wdot\vec{\omega}}
$
on the space $N=\inv{\vec{\mu}}(0)$.
\begin{claim}\label{expofnasty}
On a small enough tubular neighbourhood $U$ of $N$ we have
\begin{eqnarray*}
& &\exp{i\vec{\omega}\wdot\vec{\omega}}\\
&=&\exp{i{\vec{\omega}'}\wdot{\vec{\omega}'}
-i\vec{\omega}'\wdot\left(\phi\wedge\ol{\phi}\right)
+6i\phi_0\wedge\phi_1\wedge\phi_2\wedge \phi_3}\\
&=&\exp{i{\vec{\omega}'}\wdot{\vec{\omega}'}}\left(
1-i\vec{\omega}'\wdot\left(\phi\wedge\ol{\phi}\right)
-24{\vec{\omega}'}\wdot{\vec{\omega}'}\wedge\phi_0\wedge\phi_1\wedge\phi_2\wedge
\phi_3 +6i\phi_0\wedge\phi_1\wedge\phi_2\wedge \phi_3\right).
\end{eqnarray*}
\end{claim}
This is just a calculation.
%\begin{prf}
%Let us simplify matters by temporarily writing
%\begin{eqnarray*}
%A&=&i\vec{\omega}'\wdot\vec{\omega}'\\
%B&=&-i\vec{\omega}'\wdot\left(\phi\wedge\ol{\phi}\right)\\
%C&=&6i\phi_0\wedge\phi_1\wedge\phi_2\wedge \phi_3.
%\end{eqnarray*}
%The quantities $A,B$ and $C$ commute (they are even degree
%forms), and it is not hard to prove that satisfy
%\begin{eqnarray*}
%B^2&=&4AC,\\
%C^2&=&0,\\
% BC&=&0.
%\end{eqnarray*}
%Then
%\begin{eqnarray*}
%\left(A+B+C\right)^n
%&=&\sum_{a=0}^{n}\sum_{b=0}^a\frac{n!}{a!(n-a)!}\frac{a!}{b!(a-b)!}A^{n-a}B^{a-b}C^b\\
%&=&\sum_{a=0}^{n}\sum_{b=0}^a\frac{n!}{(n-a)!b!(a-b)!}A^{n-a}B^{a-b}C^b\\
%&=&\sum_{a=0}^{n}\frac{n!}{(n-a)!a!}A^{n-a}B^{a}+\sum_{a=1}^{n}\frac{n!}{(n-a)!(a-1)!}A^{n-a}B^{a-1}C\\
%&=&\sum_{a=0}^{2}\frac{n!}{(n-a)!a!}A^{n-a}B^{a}+nA^{n-1}C\\
%&=&A^n+nA^{n-1}B+\frac{1}{2}n(n-1)A^{n-2}B^2+nA^{n-1}C.\\
%%&=&A^n+nA^{n-1}B+4n(n-1)A^{n-2}B^2+nA^{n-1}C.
%\end{eqnarray*}
%Hence
%\begin{eqnarray*}
%\exp{A+B+C}&=&1+\sum_{n=1}^\infty\frac{\left(A+B+C\right)^n}{n!}\\
%           &=&1+\sum_{n=1}^\infty\frac{A^n+nA^{n-1}B+\frac{1}{2}n(n-1)A^{n-2}B^2+nA^{n-1}C.}{n!}\\
%&=&1
%+\sum_{n=1}^\infty\frac{A^n}{n!}
%+\sum_{n=1}^\infty\frac{A^{n-1}B}{(n-1)!}
%+\sum_{n=2}^\infty\frac{A^{n-2}B^2}{(n-2)!}
%+\sum_{n=1}^\infty\frac{A^{n-1}C.}{(n-1)!}\\
%&=&\exp{A}+\exp{A}B+4\exp{A}AC+\exp{A}C\\
%&=&\exp{A}\left(1+B+4AC+C\right),
%\end{eqnarray*}
%and substituting back gives the desired result.
%\end{prf}

Since the normal bundle of $N$ is trivialised by the nowhere vanishing sections
$\vec{q}X_\xi$ for $\vec{q}\in\{\I,\J,\K\}$, we can write $U=N\times V_\eps$
where $V_\eps$ is the ball in $\Im\quat$ centred at $0$ with radius $\epsilon$ and such that if $x=(\hat{x},\vec{z})\in N\times V_\eps$, then in this trivialisation $\vec{\mu}(\hat{x},\vec{z})=\vec{z}$. This follows from the fact that for
any local submersion $f$ of manifolds, there are local coordinate charts such
that $f$ is locally a projection, see \cite{GP}.
So we if we set
\[
W(\vec{z},\sqrt{y})=\exp{2i\sqrt{y}\vec{z}.\vec{\omega}},
\]
we can write
\begin{eqnarray*}
& &
\oint_{M}\exp{i\vec{\omega}\wdot\vec{\omega}+2i{\sqrt{y}}{\vec{\mu}}.\vec{\omega}+i|{\vec{\mu}}|^2y}\eta(\sqrt{y})\\
&=&\oint_{N\times V_\eps}\exp{i\vec{\omega}'\wdot\vec{\omega}'-i\vec{\omega}'\wdot\left(\phi\wedge\ol{\phi}\right)+6i\phi_0\wedge\phi_1\wedge\phi_2\wedge \phi_3+i|{\vec{z}|^2y}}W(\vec{z},\sqrt{y})\iota^*\eta(\sqrt{y})\\
&=&\oint_{N\times V_\eps}\exp{i\vec{\omega}'\wdot\vec{\omega}'-i\vec{\omega}'\wdot\left(\phi\wedge\ol{\phi}\right)+6i\phi_0\wedge\phi_1\wedge\phi_2\wedge \phi_3+i|{\vec{z}|^2y}}W(\vec{z},\sqrt{y})\left(\pi^*\eta_0+\dg\beta(\sqrt{y})\right)\\
&=&\oint_{N\times V_\eps}\exp{i\vec{\omega}'\wdot\vec{\omega}'-i\vec{\omega}'\wdot\left(\phi\wedge\ol{\phi}\right)+6i\phi_0\wedge\phi_1\wedge\phi_2\wedge \phi_3+i|{\vec{z}|^2y}}W(\vec{z},\sqrt{y})\pi^*\eta_0\\
\end{eqnarray*}
using a similar argument in Theorem \ref{genequivintegral} and the
compactness of the support of $\beta$. Using the fact that
$\d\vec{z}=X_{i}\hook\vec{\omega}=\phi-\phi_0$, we see that
\begin{eqnarray*}
& &\oint_{M}\exp{i\vec{\omega}\wdot\vec{\omega}+2i{\sqrt{y}}{\vec{\mu}}.\vec{\omega}+i|{\vec{\mu}}|^2y}\eta(\sqrt{y})\\
&=&\oint_{N\times V_\eps}\exp{i\vec{\omega}'\wdot\vec{\omega}'-i\vec{\omega}'\wdot\left({\phi}\wedge\ol{\phi}\right)+6i\phi_0\wedge\d vol(\vec{z})+i|{\vec{z}|^2y}}W(\vec{z},\sqrt{y})\pi^*\eta_0\\
&=&\oint_{N\times V_\eps}\exp{i\vec{\omega}'\wdot\vec{\omega}'+i|{\vec{z}|^2y}}W(\vec{z},\sqrt{y})\pi^*\eta_0\\
&-&i\oint_{N\times V_\eps}\exp{i\vec{\omega}'\wdot\vec{\omega}'+i|{\vec{z}|^2y}}\vec{\omega}'\wdot\left(\d\vec{z}\wedge{\d{\vec{z}\,}^*}\right)\wedge W(\vec{z},\sqrt{y})\pi^*\eta_0\\
&-&24\oint_{N\times V_\eps}\exp{i\vec{\omega}'\wdot\vec{\omega}'+i|{\vec{z}|^2y}}\vec{\omega}'\wdot\vec{\omega}'\wedge \phi_0\wedge\d vol(\vec{z})\wedge W(\vec{z},\sqrt{y})\pi^*\eta_0\\
&+&6i\oint_{N\times V_\eps}\exp{i\vec{\omega}'\wdot\vec{\omega}'+i|{\vec{z}|^2y}}\phi_0\wedge\d vol(\vec{z})\wedge W(\vec{z},\sqrt{y})\pi^*\eta_0\\
&=&6i\oint_{N\times V_\eps}\exp{i\vec{\omega}'\wdot\vec{\omega}'+i|{\vec{z}|^2y}}\left(4i\vec{\omega}'\wdot\vec{\omega}+1\right)\phi_0\wedge\d vol(\vec{z})\wedge W(\vec{z},\sqrt{y})\pi^*\eta_0\\
\end{eqnarray*}
where we have also used the Claim \ref{expofnasty}.
Now we examine the integrals in $\vec{z}$ and $y$.
The integral has the form
\begin{eqnarray*}
\lim_{t\mapping\infty}\int_{V_\eps}\int_\reals\exp{-\frac{y^2}{4t}+i|{\vec{z}|^2y}}W(\vec{z},\sqrt{y})\d
vol(\vec{z})\d y.
\end{eqnarray*}
Let $\Sigma_+$ be the upper half unit-hemisphere in $\Im\quat$ and write
\[
\vec{z}=r\vec{\theta}
\]
where $\vec{\theta}\in\Sigma_+$.
\begin{eqnarray*}
& &
\lim_{t\mapping\infty}\int_{V_\eps}\int_\reals\exp{-\frac{y^2}{4t}+i|{\vec{z}|^2y}}W(\vec{z},\sqrt{y})\d
vol(\vec{z})\d y\\
&=&\int_{-\eps}^{\eps}\int_{\Sigma_+}\int_{-\infty}^\infty\exp{ir^2y}W(r\vec{\theta},\sqrt{y})r^2\d
y\d vol(\vec{\theta})\d r.
\end{eqnarray*}
Now
\[
W(\vec{z},\sqrt{y})=\exp{2i\sqrt{y}\vec{z}.\vec{\omega}},
\]
and hence is a polynomial in a term of the form $\sqrt{y}\vec{z}.\vec{a}$ for fixed quaternion $\vec{a}$.
It makes sense then to evaluate the integral
\[
\int_{-\eps}^{\eps}\int_{\Sigma_+}\int_{-\infty}^\infty\exp{ir^2y}(\sqrt{y}\vec{\theta}.\vec{a})^nr^{n+2}\d y\d vol(\vec{\theta})\d r
\]
for $n\in\naturals$.
Notice that when $n$ is odd, the integral necessarily vanishes as the integral
in $r$ is the integral of an odd function. This leaves us with the case $n$
even.
So
\begin{eqnarray*}
& &
\int_{-\eps}^{\eps}\int_{\Sigma_+}\int_{-\infty}^\infty\exp{ir^2y}(\sqrt{y}\vec{\theta}.\vec{a})^2nr^{2n+2}\d
y\d vol(\vec{\theta})\d r\\
&=&\int_{-\eps}^{\eps}\int_{\Sigma_+}\int_{-\infty}^\infty\exp{ir^2y}
y^n
(\vec{\theta}.\vec{a})^{2n} r^{2n+2}\d y\d vol(\vec{\theta})\d r\\
&=&\left(\int_{\Sigma_+}(\vec{\theta}.\vec{a})^{2n}\d vol(\vec{\theta})\right)
\int_{-\eps}^{\eps}\int_{-\infty}^\infty\exp{ir^2y} y^nr^{2n+2}\d y\d r\\
&=& const\int_{-\eps}^{\eps}r^{2n+2}D^n(\delta)(r)\d r\\
&=& D^n(r\mapsto r^{2n+2})(0)\\
&=& 0.
\end{eqnarray*}
Thus
\begin{eqnarray*}
& &
\lim_{t\mapping\infty}\int_{V_\eps}\int_\reals\exp{-\frac{y^2}{4t}+i|{\vec{z}|^2y}}W(\vec{z},\sqrt{y})\d
vol(\vec{z})\d y\\
&=&
\lim_{t\mapping\infty}\int_{V_\eps}\int_\reals\exp{-\frac{y^2}{4t}+i|{\vec{z}|^2y}}W(\vec{z},0)\d vol(\vec{z})\d y\\
&=&\lim_{t\mapping\infty}\int_{V_\eps}\int_\reals\exp{-\frac{y^2}{4t}+i|{\vec{z}|^2y}}\d vol(\vec{z})\d y\\
&=&2\pi\int_{V_\eps}\delta(|\vec{z}|^2)\d vol(\vec{z})\\
&=&2\pi
\end{eqnarray*}
and
\begin{eqnarray*}
& &\oint_M
\exp{i\vec{\omega}\wdot\vec{\omega}+2i{\sqrt{y}}{\vec{\mu}}.\vec{\omega}+i
|{\vec{\mu}}|^2y}\eta(\sqrt{y})\\
&=&6i\frac{2\pi}{vol\sphere{1}}\oint_{N}\exp{i\vec{\omega}'\wdot\vec{\omega}'}\left(4i\vec{\omega}'\wdot\vec{\omega}+1\right)\phi_0\wedge\pi^*\eta_0\\
&=&6i\frac{2\pi}{vol\sphere{1}}\oint_{N}\pi^*\left(\exp{i\vec{\omega_0}\wdot\vec{\omega_0}}\left(4i\vec{\omega_0}\wdot\vec{\omega_0}+1\right)\right)\phi_0\wedge\pi^*\eta_0.\\
\end{eqnarray*}
Now, recall that $\phi_0$ is $\frac{1}{\sqrt{2}}$ times the connection 1-form which acts as a fibre volume form.
Thus we have
\begin{eqnarray*}
& &\oint_M
\exp{i\vec{\omega}\wdot\vec{\omega}+2i{\sqrt{y}}{\vec{\mu}}.\vec{\omega}+i
|{\vec{\mu}}|^2y}\eta(\sqrt{y})\\
&=&6i\frac{2\pi}{\sqrt{2}}\oint_{\mathcal{M}}\exp{i\vec{\omega_0}\wdot\vec{\omega_0}}\left(4i\vec{\omega_0}\wdot\vec{\omega_0}+1\right)\wedge\eta_0\\
&=&6i\pi\sqrt{2}\oint_{\mathcal{M}}\exp{i\vec{\omega_0}\wdot\vec{\omega_0}}
\left(4i\vec{\omega_0}\wdot\vec{\omega_0}+1\right)\wedge\eta_0.\\
&=&6i\pi\sqrt{2}(\dim\mathcal{M}-\mathrm{deg}\eta_0+1)\oint_{\mathcal{M}}\exp{i\vec{\omega_0}\wdot\vec{\omega_0}}
\wedge\eta_0.\\
\end{eqnarray*}
If the right hand side is known to exist, then the above argument
may be reversed.
\end{prf}
\begin{cor1}\label{hypintbound}
If $M$ is a hyperK\"ahler manifold with a tri-Hamiltonian action
of $\sphere{1}$ and $L$ an invariant submanifold respected by
$\sphere{1}$ . Further, if
$$\mathcal{L}=(L\cap\vec{\mu}^{-1}(0))/\sphere{1}$$ and
$\eta\in\Omega_{\sphere{1}}^\bullet(L)$
is compactly supported and compactly associated with
$\eta_0\in\FORM{\mathcal{L}}{\bullet}$, then
\begin{eqnarray*}
& &
(\dim\mathcal{L}-\mathrm{deg}\eta_0+1)\int_{\mathcal{L}}\exp{i\iota_0^*(\vec{\omega}_0\wdot\vec{\omega}_0)}
\eta_0\\
&=&\frac{1}{6\pi i\sqrt{2}}\oint_{L}
\exp{i\iota^*(\vec{\omega}\wdot\vec{\omega}+2i{\sqrt{y}}{\vec{\mu}}.\vec{\omega}+i
|{\vec{\mu}}|^2y)}\eta(\sqrt{y}),
\end{eqnarray*}
where $\iota:L \hookrightarrow M$ and
$\iota_0:\mathcal{L}\hookrightarrow\mathcal{M}$ are the inclusions
.
\end{cor1}
\begin{prf}
This follows the proof of the above theorem, the main consequence
of the action respecting the submanifold is that the Poincar\'e
dual of $\vec{\mu}^{-1}(0)$ in $M$ does not vanish when we
restrict it to $L$.
\end{prf}

\begin{rmk}
In these integrals, we have not assumed anything about the (equivariant) closure of $\eta$.
The integrals work well enough without that assumption. However, in order to
use localisation, we \emph{do} need assumptions on equivariant closure.
\end{rmk}
\subsection{Localisation of HyperK\"ahler Integrals for $G=\sphere{1}$}
At first sight, it appears that a localisation formula would be difficult to
apply to the Witten-style equation we derived in Theorem \ref{hypintform} due
to the presence of the square roots. What is most important in this formula is
that it doesn't depend on choice of square root! We can exchange $-\sqrt{y}$
for $\sqrt{y}$ in the formula and obtain the same result%
\footnote{Notice that we cannot replace $-y$ for $y$ in the symplectic formula
without changing the result. This is due to the fact that by changing the sign
in this situation, we change a Fourier transform into an inverse Fourier transform}.
Hence we have the result for the circle
\begin{thm1}\label{hypintform2}
\begin{eqnarray*}
& &
(\dim\mathcal{M}-\mathrm{deg}\eta_0+1)\int_\mathcal{M}\exp{i\vec{\omega}_0\wdot\vec{\omega}_0}
\eta_0\\
&=&\frac{1}{12i\pi\sqrt{2}}\oint_M
\exp{i\vec{\omega}\wdot\vec{\omega}+i|{\vec{\mu}}|^2y}\left(\exp{2i{\sqrt{y}}{\vec{\mu}}.\vec{\omega}}\eta(\sqrt{y})+\exp{-2i{\sqrt{y}}{\vec{\mu}}.\vec{\omega}}\eta(-\sqrt{y})\right).
\end{eqnarray*}
\end{thm1}
\begin{prf}
Follows from Theorem \ref{hypintform} by summing
\begin{eqnarray*}
(\dim\mathcal{M}-\mathrm{deg}\eta_0+1)\int_\mathcal{M}\exp{i\vec{\omega}_0\wdot\vec{\omega}_0}\eta_0
&=&\frac{1}{6\pi i\sqrt{2}}\oint_M
\exp{i\vec{\omega}\wdot\vec{\omega}+2i{\sqrt{y}}{\vec{\mu}}.\vec{\omega}+i
|{\vec{\mu}}|^2y}\eta(\sqrt{y})
\end{eqnarray*}
and
\begin{eqnarray*}
(\dim\mathcal{M}-\mathrm{deg}\eta_0+1)\int_\mathcal{M}\exp{i\vec{\omega}_0\wdot\vec{\omega}_0}\eta_0
&=&\frac{1}{6\pi i\sqrt{2}}\oint_M
\exp{i\vec{\omega}\wdot\vec{\omega}-2i{\sqrt{y}}{\vec{\mu}}.\vec{\omega}+i
|{\vec{\mu}}|^2y}\eta(-\sqrt{y})
\end{eqnarray*}
to get
\begin{eqnarray*}
& &
2(\dim\mathcal{M}-\mathrm{deg}\eta_0+1)\int_\mathcal{M}\exp{i\vec{\omega}_0\wdot\vec{\omega}_0}\eta_0\\
&=&\frac{1}{6\pi i\sqrt{2}}\oint_M
\exp{i\vec{\omega}\wdot\vec{\omega}+i|{\vec{\mu}}|^2y}\left(\exp{2i{\sqrt{y}}{\vec{\mu}}.\vec{\omega}}\eta(\sqrt{y})+\exp{-2i{\sqrt{y}}{\vec{\mu}}.\vec{\omega}}\eta(-\sqrt{y})\right).
\end{eqnarray*}
\end{prf}
Now for any polynomial $P\in\complexes[X]$, we know that the polynomial
\[
Q:y\mapsto P(y)+P(-y)
\]
is a polynomial consisting only of even powers. Hence $Q(\sqrt{y})$ is a
perfectly defined polynomial in $y$. By defining
\begin{eqnarray*}
&\mathfrak{P}:\complexes[X]\mapping\complexes[X^2]\\
&\mathfrak{P}(Q)(X)=\frac{1}{2}\left(Q(X)+Q(-X)\right)
\end{eqnarray*}
we see that
\begin{equation}\label{circwitt}
(\dim\mathcal{M}-\mathrm{deg}\eta_0+1)\int_\mathcal{M}\exp{i\vec{\omega}_0\wdot\vec{\omega}_0}\eta_0=\frac{1}{6\pi
i\sqrt{2}}\oint_M
\exp{i\vec{\omega}\wdot\vec{\omega}+i|{\vec{\mu}}|^2y}\mathfrak{P}\left(z\mapsto\exp{2i{z}{\vec{\mu}}.\vec{\omega}}\eta(z)\right)(\sqrt{y}).
\end{equation}
Examining the right hand side of the formula
in Theorem \ref{hypintform2}, we can see that the integrand in
Theorem \ref{hypintform2} is an entire function of $y$. This
allows us to apply the theory that Jeffrey and Kirwan develop in
\cite{JK} to reduce the study after localising to residue
formul\ae. After applying this we find that in the case for the
circle we have the following localisation formula using a similar
proof to that of Lemma \ref{Resthm0}.
\begin{thm1}\label{circlocform},
Suppose $\eta\in\Omega_{\sphere{1}}^\bullet(M)$ is equivariantly closed and associated with
the closed form $\eta_0\in\FORM{\mathcal{M}}{\bullet}$.
 Let $\iota:M_0\hookrightarrow M$ be the inclusion of the fixed point set
 in $M$ and $e$ be the equivariant Euler class of the normal bundle of $M_0$ in $M$.
Then
\begin{eqnarray*}
& &(\dim\mathcal{M}-\mathrm{deg}\eta_0+1)\int_\mathcal{M}\exp{i\vec{\omega}_0\wdot\vec{\omega}_0}\eta_0\\
&=&\frac{-1}{6\sqrt{2}\pi\mathrm{vol}\,\sphere{1}}\int_{M_0}\mathrm{Coeff}_{y^{-1}}
\left[\iota^*\exp{i\vec{\omega}\wdot\vec{\omega}+i|{\vec{\mu}}|^2y}\mathfrak{P}\left(z\mapsto\frac{\iota^*\exp{2i{z}{\vec{\mu}}.\vec{\omega}}\eta(z)}{e(z)}\right)(\sqrt{y})\right]\\
&=&\frac{-1}{6\sqrt{2}\pi\mathrm{vol}\,\sphere{1}}\int_{M_0}\mathrm{Coeff}_{y^{-2}}
\left[\iota^*\exp{i\vec{\omega}\wdot\vec{\omega}+i|{\vec{\mu}}|^2y^2}\left(\frac{\iota^*\exp{2i{y}{\vec{\mu}}.\vec{\omega}}\eta(y)}{e(y)}\right)\right].
\end{eqnarray*}
\end{thm1}
\begin{prf}
First we examine the right hand side of (\ref{circwitt}). By Theorem
\ref{locform}, we know that
\[
\int_M
\exp{i\vec{\omega}\wdot\vec{\omega}+i|{\vec{\mu}}|^2X^2}
\exp{2i{X}{\vec{\mu}}.\vec{\omega}}\eta(X)
=\int_{M_0}\frac{\iota^*\left[
\exp{i\vec{\omega}\wdot\vec{\omega}+i|{\vec{\mu}}|^2X^2}
\exp{2i{X}{\vec{\mu}}.\vec{\omega}}\eta(X)
\right]
}{e(X)}.
\]
So we have
\begin{eqnarray*}
& &\int_M
\exp{i\vec{\omega}\wdot\vec{\omega}+i|{\vec{\mu}}|^2X^2}
\mathfrak{P}\left(z\mapsto
\exp{2i{z}{\vec{\mu}}.\vec{\omega}}\eta(z)
\right)(X)\\
&=&\int_{M_0}\iota^*\left[
\exp{i\vec{\omega}\wdot\vec{\omega}+i|{\vec{\mu}}|^2X^2}
\mathfrak{P}\left(z\mapsto
\frac{
\exp{2i{z}{\vec{\mu}}.\vec{\omega}}\eta(z)}
{e(z)}
\right)(X)
\right]
\end{eqnarray*}
hence
\begin{eqnarray*}
& &\int_\reals[\d y]\exp{-\frac{y^2}{4t}}\int_M
\exp{i\vec{\omega}\wdot\vec{\omega}+i|{\vec{\mu}}|^2y}\mathfrak{P}\left(z\mapsto\exp{2i{z}{\vec{\mu}}.\vec{\omega}}\eta(z)\right)(\sqrt{y})\\
&=&\int_\reals[\d y]\exp{-\frac{y^2}{4t}}\int_{M_0}\iota^*\left[
\exp{i\vec{\omega}\wdot\vec{\omega}+i|{\vec{\mu}}|^2y}
\mathfrak{P}\left(z\mapsto
\frac{
\exp{2i{z}{\vec{\mu}}.\vec{\omega}}\eta(z)}
{e(z)}
\right)(\sqrt{y})
\right]
\end{eqnarray*}
by putting $X=\sqrt{y}$.
 Since the integrand is a rational function of $y$, the
rest now follows from the same proof as Lemma \ref{Resthm0}
\end{prf}
\begin{rmk1}
Notice that in the symplectic case, the only contribution from the fixed
point set came from those components on which $\mu$ was positive. Here the
r\^ole of $\mu$ is taken by $|\vec{\mu}|^2$ which is always positive. It
follows that each component of $M_0$ will make a contribution in the hyperK\"ahler
case.
\end{rmk1}
Using this we can just extend this idea to a product of circles to get
a Witten-style expression and residue theorem for the torus.
\begin{thm1}\label{torgenhypint}
Let $M$ be a hyperK\"ahler manifold with tri-Hamiltonian action of
the torus $\mathbb{T}^k$. Let
$\mathcal{M}=M/\!\!/\!\!/\!\!/\mathbb{T}^k$ and with hyperK\"ahler
form $\vec{\omega}_0$, then for each compactly supported
$\eta_0\in\FORM{\mathcal{M}}{\bullet}$ compactly associated with
$\eta\in\OG{M}{\bullet}$,
we have
\begin{eqnarray*}
& &
(\dim\mathcal{M}-\mathrm{deg}\eta_0+1)\int_\mathcal{M}\exp{i\vec{\omega}_0\wdot\vec{\omega}_0}\eta_0\\
&=&\left(\frac{1}{6\pi i\sqrt{2}}\right)^k\oint_M
\exp{i\vec{\omega}\wdot\vec{\omega}+i\sum_{n=0}^k|{\vec{\mu}_\nu}|^2y_\nu}\mathfrak{P}\left(\ul{z}\mapsto\exp{2i\sum_{\nu=0}{k}{z}{\vec{\mu}_\nu}.\vec{\omega}}\eta(\ul{z})\right)(\sqrt{\ul{y}})\\
\end{eqnarray*}
where $\sqrt{\ul{y}}=(\sqrt{y_1},\ldots,\sqrt{y_k})$.
\end{thm1}
\begin{prf}
We decompose $\mathbb{T}^k$ into the product of circles, and recursively use Theorem
\ref{hypintform}.
\end{prf}
\begin{cor1}\label{hyptorloc}
\begin{eqnarray*}
& &(\dim\mathcal{M}-\mathrm{deg}\eta_0+1)\left({-6\pi\sqrt{2}\mathrm{vol}\,\sphere{1}}\right)^k\int_\mathcal{M}\exp{i\vec{\omega}_0\wdot\vec{\omega}_0}\eta_0\\
&=&\int_{M_0}\mathrm{Coeff}_{{y_1}^{-1}\ldots{y_k}^{-1}}
\left[\iota^*\exp{i\vec{\omega}\wdot\vec{\omega}+i\sum_{\nu=1}^k|{\vec{\mu}_\nu}|^2y_\nu}\mathfrak{P}\left(\ul{z}\mapsto\frac{\iota^*\exp{2i\sum_{\nu=1}^k{z_\nu}{\vec{\mu}_\nu}.\vec{\omega}}\eta(\ul{z})}{e(\ul{z})}\right)(\sqrt{\ul{y}})\right]\\
&=&\int_{M_0}\mathrm{Coeff}_{{y_1}^{-2}\ldots{y_k}^{-2}}
\left[\iota^*\exp{i\vec{\omega}\wdot\vec{\omega}+i\sum_{\nu=1}^k|{\vec{\mu}_\nu}|^2y_\nu^2}\frac{\iota^*\exp{2i\sum_{\nu=1}^k{y_\nu}{\vec{\mu}_nu}.\vec{\omega}}\eta(\ul{y})}{e(\ul{y})}\right]\\
\end{eqnarray*}
\end{cor1}
\begin{prf}
We decompose $\mathbb{T}^k$ into the product of circles, and repeat Theorem \ref{circlocform}.
\end{prf}
Now we'd really like to extend this result to more complicated Lie groups
than tori, in order to do this we must follow the strategy of Shaun
Martin \cite{MN} and relate the integrals of a hyperK\"ahler quotient by $G$
to integrals of a hyperK\"ahler quotient by its maximal torus $T$.\par
\subsection{HyperK\"ahler Quotients and Tori}
We closely follow Shaun Martin here.
Let $G$ be a compact, connected Lie group with maximal torus $T$. Let $\Delta$
be the set of roots of $T$ with $\Delta^\pm$ being the set of $\pm$ve roots.
We have the orthogonal projection $p:\g\mapping\goth{t}$, and we denote
the moment map $\vec{\mu}:M\mapping\g\tensor\Im\quat$ by $\vec{\mu}_G$
and $(p\tensor\id)\circ\vec{\mu}:M\mapping\goth{t}\tensor\Im\quat$ by $\vec{\mu}_T$. Also set
\begin{eqnarray*}
N_G&=&\inv{{\vec{\mu}}_G}(0),\\
N_T&=&\inv{{\vec{\mu}}_T}(0),\\
\mathcal{M}_G&=&N_G/G,\\
\mathcal{M}_T&=&N_T/T.
\end{eqnarray*}
We suppose that $0$ is regular for both $\vec{\mu}_G$ and $\vec{\mu}_T$.\par
Given $\alpha\in\Delta$ we set $V_\alpha$ to be its associated root space
(isomorphic to $\complexes$) and define the vector bundle
\[
L_\alpha={N_T}\times_TV_\alpha,
\]
over $\mathcal{M}_T$.
Define
\[
V_\pm=\bigoplus_{\alpha\in\Delta^\pm}L_\alpha
\]
and
\[
V=V_+\oplus V_-.
\]
We also have
\[
\begin{array}{rccl}
   &      & \iota & \\
   &N_G/T &{\hookrightarrow}&\mathcal{M}_T\\
q  &\downarrow&  & \\
   & \mathcal{M}_G &
\end{array}.
\]
\begin{prop1}\label{adaptedSM}[cf Proposition 1.2 of \cite{MN}]
\begin{enumerate}
\item The vector bundle $V_+\mapping\mathcal{M}_T$ has a section $s$, which
is transverse to the zero section $Z$, and such that the zero set of $s$ is
$N_G/T$. It follows that the normal bundle
\[
\mathcal{V}(N_G/T;\mathcal{M}_T)\iso\iota^*V_+\tensor\Im\quat
\]
\item Let $\mathrm{Vert(q)}$ be the vertical subspace of the fibration
$q:N_G/T\mapping\mathcal{M}_G$, that is the kernel of $q_*$. Then
\[
\mathrm{Vert}(q)\iso\iota^*V_+
\]
\end{enumerate}
\end{prop1}
This is a simple extension of Martin's results in \cite{MN}.
\begin{thm1}\label{torusquot}
Let $e_+=e(V_+)$, then for each compactly supported
$\alpha\in\FORM{\mathcal{M}_G}{\bullet}$ with lift
$\tilde\alpha\in\FORM{\mathcal{M}_T}{\bullet}$ (that is
$\iota^*\tilde{\alpha}=q^*\alpha$)
\[
\int_{\mathcal{M}_G}\alpha=\frac{1}{|W|}\int_{\mathcal{M}_T}\tilde{\alpha}\wedge
e_+^4,
\]
where $W$ is the Weyl group of $T$ in $G$.
\end{thm1}
\begin{prf}(Adapted from the proof of Theorem B \cite{MN})\\
First note that $\iota^* e_+=e(\mathrm{Vert(q)})$ by \ref{adaptedSM}. By
arguments given in \cite{MN}, $q_*\iota^*e_+=|W|$. Thus
\begin{eqnarray*}
\int_{\mathcal{M}_G}\alpha&=&\frac{1}{|W|}\int_{\mathcal{M}_G} \alpha\wedge
q_*\iota^*e_+\\
&=& \frac{1}{|W|}\int_{N_G/T} q^*\alpha\wedge\iota^*e_+\\
&=& \frac{1}{|W|}\int_{N_G/T} \iota^*\tilde{\alpha}\wedge\iota^*e_+\\
&=& \frac{1}{|W|}\int_{\mathcal{M}_T} \iota_*\iota^*(\tilde\alpha\wedge e_+)\\
&=& \frac{1}{|W|}\int_{\mathcal{M}_T} \tilde{\alpha}\wedge e_+\wedge e(V_+\tensor\Im\quat)\\
&=& \frac{1}{|W|}\int_{\mathcal{M}_T}\tilde{\alpha}\wedge e_+^4.
\end{eqnarray*}
\end{prf}
\begin{lem1}
The equivariant representative of $e_+\in\mathrm{H}^\bullet(\mathcal{M}_T)$ is given by
\[
w(y)=\prod_{\alpha\in\Delta^+}\alpha(y)\in\mathrm{H}_T^\bullet(M).
\]
\end{lem1}
\begin{prf}
Since $\pi_\alpha:L_\alpha=N_T\times_T V_\alpha\mapping \mathcal{M}_T$ is a line bundle, we may form $q^*L_\alpha\mapping N_T$ where $q:N_T\mapping\mathcal{M}_T$ is the quotient map. From the standard theory of principal fibrations,
we know that
\[
q^*L_\alpha\iso N_T\times L_\alpha.
\]
This is not trivial in the equivariant sense, so we introduce the equivariant
connection $\d+\alpha$. Hence
\[
q^*e(L_\alpha)(\xi)=e(q^*L_\alpha)(\xi)=\alpha(\xi).
\]
Thus we may conclude that $e(L_\alpha)(\xi)=\alpha(\xi).$ \par
Hence $$e(V_+)=\prod_{\alpha\in\Delta_+}e(L_\alpha)=\prod_{\alpha\in\Delta_+}\alpha.$$
\end{prf}
Finally we see that
\begin{thm1}\label{HYPINTFORM}
If $\eta$ is compactly associated with $\eta_0$, then
\begin{eqnarray*}
& &(\dim\mathcal{M}-\mathrm{deg}\eta_0+1))\int_\mathcal{M}\exp{i\vec{\omega}_0\wdot\vec{\omega}_0}\eta_0\\
&=&
\left(\frac{1}{6\pi i\sqrt{2}}\right)^k\frac{1}{|W|}\oint_M
\exp{i\vec{\omega}\wdot\vec{\omega}+i|{\vec{\mu}}|^2y}\mathfrak{P}\left(z\mapsto\exp{2i{z}{\vec{\mu}}.\vec{\omega}}w(z)^4\eta(z)\right)(\sqrt{y}).\\
\end{eqnarray*}
\end{thm1}
The consequence of this theorem is that now we can replace the Lie group $G$
by its maximal torus and use \ref{hypintform} recursively \`a la Guillemin and Kalkman\cite{GK}.
\subsection{Concluding Remarks}
The techniques used here apply also for Quaternionic K\"ahler manifolds.
Theorem \ref{HNormalsympform} works due to the fact that, although the complex structures are only locally defined, the 2-form $\vec{\omega}$ with values
in $\End{\Tang{M}}$ and the 4-form $\Omega=\vec{\omega}\wdot\vec{\omega}$
exist globally with identical relations on the moment \emph{section}, $\vec{\mu}$. \par
If $S\subset M$ is a 3-Sasakian hypersurface whose structure is induced by
the hyperK\"ahler structure on $M$ which posesses a trihamiltonian group
action, then we may obtain a 3-Sasakian reduction of $S$. The integration
formul\ae\ derived here will be valid for the reduction of $S$ and will prove
useful in calculating intersection forms since it is possible to have a compatible
action on a compact 3-Sasakian manifold.
\par
The Author would like to thank the EPSRC for his financial support, the American
Institute of Mathematics for its helpful workshop on ``Moment maps and Surjectivity
in various Geometries" as well as the organisers Dr. Eugene Lerman, Dr. Tara Holm
and Dr. Susan Tolman, and is appreciative of (in
no particular order) the very welcome input of Dr. Alistair Craw, Prof. Frances Kirwan, Prof. Lisa Jeffrey, Prof. John Rawnsley, Dr. Roger Bielawski and Dr. Richard Thomas,
and Prof. Charles Boyer. HOAMGD.

\end{document}